\documentclass{article}

\usepackage{kpfonts} 
\usepackage[T1]{fontenc}

\title{Insights from a workshop on gamification of research in mathematics and computer science}
\author{Alexis Langlois-Rémillard, Élise Raphael, Erika Roldan}
\date{\today}

\AtEndDocument{\bigskip{\footnotesize%
  (A.~Langlois-R\'emillard)  \textsc{Hausdorff Center for Mathematics, Universit\"at Bonn, Bonn, Germany}\par
  \addvspace{\medskipamount}
  Email: \href{langlois@uni-bonn.de}{langlois@uni-bonn.de} website: \href{https://alexisl-r.github.io/}{https://alexisl-r.github.io/}  ORCID: \href{https://orcid.org/0000-0002-5919-8766}{0000-0002-5919-8766}\par
  \addvspace{\medskipamount}
  (\'E.~Raphael) \textsc{Universit\'e de Gen\`eve, Geneva, Switzerland} \par
    \addvspace{\medskipamount}
\href{Elise.Raphael@unige.ch}{Elise.Raphael@unige.ch}\par
  \addvspace{\medskipamount}
 (E.~Rold\'an ) \textsc{Max Planck Institute for Mathematics in the Sciences, and ScaDS.AI, Universit\"at Leipzig, Leipzig, Germany.}\par
 Email: \href{roldan@mpg.mis.de}{roldan@mpg.mis.de} website: \href{https://www.erikaroldan.net/}{https://www.erikaroldan.net/} ORCID: \href{https://orcid.org/0000-0002-4788-8225}{0000-0002-4788-8225}
} 
}

\usepackage{amsthm}

\newenvironment{dialogue}{\list{---}{\itemsep=\parskip \topsep=\parskip \parsep=\parskip}}{\endlist}

\usepackage{hyperref}
\hypersetup{
colorlinks=true,
citecolor=teal,
linkcolor=teal,
urlcolor=cyan
}
\usepackage{geometry}
\usepackage[dvipsnames]{xcolor}

\usepackage{tikz}
\usepackage{tikz-cd}
\usetikzlibrary{arrows}
\usetikzlibrary{shapes.arrows}   
\usetikzlibrary{positioning}
\usepackage{tcolorbox}

\usepackage[hyperref,backref,backrefstyle=none,giveninits,block=space]{biblatex}
\renewbibmacro{in:}{} 

\usepackage{subcaption}

\usepackage{wrapfig} 

\usepackage[capitalise,noabbrev]{cleveref}
\newcommand{\opacite}{.8}
\newcommand{\radis}{0pt}
\tikzset
{%
    pics/pentominoV/.style={code={\draw[thick,fill=\couleura,fill opacity=\opacite] (0,0) |- (1,3) |- (3,1) |- cycle;
              \fill[red] (0,0) circle [radius=\radis]; 
              }},
    pics/pentominoX/.style={code={\draw[thick,fill=\couleurb,fill opacity=\opacite] (1,0) |- (0,1) |- (1,2) |- (2,3) |- (3,2) |- (2,1) |- cycle;
    \fill[red] (0,0) circle [radius=\radis]; 
    }},
    pics/pentominoT/.style={code={\draw[thick,fill=\couleurc,fill opacity=\opacite] (1,0) |- (0,2) |- (3,3) |- (2,2) |- cycle;
    \fill[red] (0,0) circle [radius=\radis]; 
    }}, 
    pics/pentominoN/.style={code={\draw[thick,xscale=#1,
                                        fill=\couleurd ,fill opacity=\opacite] (0,0) |- (1,3) |- (2,4) |- (1,2) |- cycle;
                                  \fill[red] (0,0) circle [radius=\radis]; 
                                  }},
    pics/pentominoN/.default={1},
    pics/pentominoU/.style={code={\draw[thick,fill=\couleure ,fill opacity=\opacite] (0,0) |- (1,2) |- (2,1) |- (3,2) |- cycle;
    \fill[red] (0,0) circle [radius=\radis];
    }},
    pics/pentominoW/.style={code={\draw[thick,fill=\couleurf,fill opacity=\opacite] (0,0) |- (1,1) |- (2,2) |- (3,3) |- (2,1) |- cycle;
              \fill[red] (0,0) circle [radius=\radis]; 
              }},
   pics/pentominoZ/.style={code={\draw[xscale=#1,
                                        thick,fill=\couleurg ,fill opacity=\opacite] (0,0) |- (1,1) |- (3,3) |- (2,2) |- cycle;
                                  \fill[red] (0,0) circle [radius=\radis]; 
                                  }},
    pics/pentominoZ/.default={1},
     pics/pentominoI/.style={code={\draw[thick,fill=\couleurh,fill opacity=\opacite] (0,0) |- (1,5) |- cycle;
      \fill[red] (0,0) circle [radius=\radis]; 
     }},
    pics/pentominoF/.style={code={\draw[xscale=#1,
                                        thick,fill=\couleuri ,fill opacity=\opacite] (1,0) |- (0,1) |- (1,2) |- (3,3) |- (2,2) |- cycle;
                                  \fill[red] (0,0) circle [radius=\radis]; 
                                  }},
    pics/pentominoF/.default={1},
    pics/pentominoY/.style={code={\draw[xscale=#1,
                                        thick,fill=\couleurj ,fill opacity=\opacite] (0,0) |- (1,4) |- (2,3) |- (1,2) |- cycle;
                                  \fill[red] (0,0) circle [radius=\radis]; 
                                 }},
    pics/pentominoY/.default={1},
    pics/pentominoP/.style={code={\draw[xscale=#1,
                                        thick,fill=\couleurk ,fill opacity=\opacite] (0,0) |- (2,3) |- (1,1) |- cycle;
                                  \fill[red] (0,0) circle [radius=\radis]; 
                                  }},
    pics/pentominoP/.default={1},
    pics/pentominoL/.style={code={\draw[xscale=#1,
                                        thick,fill=\couleurl ,fill opacity=\opacite] (0,0) |- (1,4) |- (2,1) |-  cycle;
                                  \fill[red] (0,0) circle [radius=\radis]; 
                                  }},
    pics/pentominoL/.default={1},
    pics/tetrominoL/.style={code={\draw[xscale=#1,
                                        thick,fill=\couleurl ,fill opacity=\opacite] (0,0) |- (1,3) |- (2,1) |-  cycle;
                                  \fill[red] (0,0) circle [radius=\radis]; 
                                  }},
    pics/tetrominoL/.default={1},
    pics/tetrominoI/.style={code={\draw[thick,fill=\couleurh,fill opacity=\opacite] (0,0) |- (1,4) |- cycle;
      \fill[red] (0,0) circle [radius=\radis]; 
     }},
        pics/tetrominoN/.style={code={\draw[xscale=#1,
                                        thick,fill=\couleurj ,fill opacity=\opacite] (0,0) |- (1,1) |- (3,2) |- (2,1) |- cycle;
                                  \fill[red] (0,0) circle [radius=\radis]; 
                                  }},
    pics/tetrominoN/.default={1},
    pics/tetrominoT/.style={code={\draw[thick,fill=\couleurc,fill opacity=\opacite] (1,0) |- (0,1) |- (3,2) |- (2,1) |- cycle;
    \fill[red] (0,0) circle [radius=\radis]; 
    }}, 
    pics/tetrominoO/.style={code={\draw[thick,fill=\couleure,fill opacity=\opacite] (0,0) |- (2,2)  |- cycle;
    \fill[red] (0,0) circle [radius=\radis]; 
    }}, 
   }

\newcommand{\couleura}{Goldenrod}
\newcommand{\couleurb}{JungleGreen}
\newcommand{\couleurc}{CadetBlue}
\newcommand{\couleurd}{Turquoise}
\newcommand{\couleure}{Apricot}
\newcommand{\couleurf}{DarkOrchid}
\newcommand{\couleurg}{Sepia}
\newcommand{\couleurh}{Gray}
\newcommand{\couleuri}{RawSienna}
\newcommand{\couleurj}{OliveGreen}
\newcommand{\couleurk}{BrickRed}
\newcommand{\couleurl}{WildStrawberry}

\newcommand{\scfa}{.4} 
\newcommand{\polfont}[1]{\texttt{#1}}

\theoremstyle{plain}
\newtheorem{theorem}{Theorem}[section]
\newtheorem{lemma}[theorem]{Lemma}
\newtheorem{proposition}[theorem]{Proposition}

\newtheorem{definition}[theorem]{Definition}

\usepackage{setspace} 
\begin{document}


\maketitle

\begin{abstract}
Can outreach inspire and lead to research and vice versa? In this work, we introduce our approach to the \textit{gamification of research} in mathematics and computer science through three illustrative examples. We discuss our primary motivations and provide insights into what makes our proposed gamification effective for three research topics in discrete and computational geometry and topology: \textbf{(1) DominatriX}, an art gallery problem involving polyominoes with rooks and queens; \textbf{(2) Cubical Sliding Puzzles}, an exploration of the discrete configuration spaces of sliding puzzles on the $d$-cube with topological obstructions; and \textbf{(3) The Fence Challenge}, a participatory isoperimetric problem based on polyforms. Additionally, we report on the collaborative development of the game \textit{Le Carré du Diable}, inspired by The Fence Challenge and created during the workshop \textit{Let's talk about outreach!}, held in October 2022 in Les Diablerets, Switzerland. All of our outreach encounters and creations are designed and curated with an inclusive culture and a strong commitment to welcoming the most diverse audience possible.

\end{abstract}

\tableofcontents

\section{Introduction}

To set the stage, we begin this contribution by introducing one of our research topics through a small puzzle. After all, if you are reading this article, there is a good chance you enjoy learning about and solving puzzles.

\begin{tcolorbox}[colback=gray!5!white,colframe=red!75!black,title= The Fence Challenge with tetrominoes]
Using the following five pieces (called tetrominoes), enclose the biggest area possible.
\begin{equation}
\begin{tikzpicture}[scale = \scfa,baseline={(current bounding box.center)}]
    \draw[help lines] (0,0) grid (18,6);
        \pic[scale=\scfa] at (1,1) {tetrominoI};
        \draw (1.5,1) node[below] {\polfont i};
        \pic[scale=\scfa] at (4,1) {tetrominoL};
        \draw (4.5,1) node[below] {\polfont l};
        \pic[scale=\scfa] at (7,1) {tetrominoN};
        \draw (8.5,1) node[below] {\polfont n};
        \pic[scale=\scfa] at (11,1) {tetrominoO};
        \draw (12.5,1) node[below] {\polfont o};
          \pic[scale=\scfa] at (14,1) {tetrominoT};
        \draw (15.5,1) node[below] {\polfont t};
\end{tikzpicture}
\label{eq:tetrapieces}
\end{equation}
\end{tcolorbox}
Using the five pieces depicted above in \eqref{eq:tetrapieces}, which you might recognise from the game \textit{Tetris} \cite{pajitnov1984tetris}, try to enclose the largest possible area on a square grid. To clarify the objective of the puzzle, we present a hypothetical dialogue between an imaginary reader and us.

\begin{wrapfigure}{L}{0.35\textwidth}
\centering
    \begin{tikzpicture}[scale = \scfa]
    \pic[scale=\scfa] at (-1,0) {tetrominoO};
    \pic[scale=\scfa] at (1,2) {tetrominoN};
    \draw (5,0) -- (5,4);
    \pic[scale=\scfa] at (7,0) {tetrominoO};
    \pic[scale=\scfa] at (8,2) {tetrominoN};
    \draw[very thick, ->] (7.5,1.5) -- (7.5,0.5) -- (8.5,0.5) -- (8.5,1.5) -- (8.5,2.5) --  (9.5,2.5) -- (9.5,3.5) -- (10.5,3.5);
\end{tikzpicture}\caption{Left, a corner-connection, not closing the fence; right an edge-connection with the rook-movement}\label{fig:corner_connect}
\end{wrapfigure}
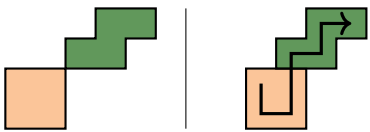
\begin{wrapfigure}{R}{0.35\textwidth}
    \centering
    \begin{tikzpicture}[scale=\scfa]
    \draw[help lines] (0,0) grid (8,8);
        \pic[scale = \scfa,rotate=270] at (1,1){tetrominoL={-1}};
        \pic[scale=\scfa,rotate=0] at (1,3) {tetrominoI};
        \pic[scale = \scfa,rotate =0] at (5,5) {tetrominoN={-1}};
        \pic[scale = \scfa,rotate=270] at (2,5){tetrominoT};
    \end{tikzpicture}
    \caption{A tentative of tetromino fence (without the ``\texttt{o}'').}
    \label{fig:tentative_tetra}
\end{wrapfigure}
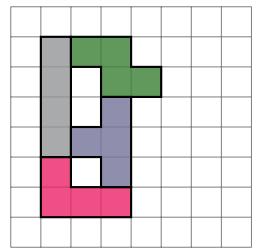

\begin{dialogue}
    \item \textbf{Reader:} The puzzle asks us to enclose the largest area, but how do we define ``enclosing''?
    \item \textbf{Us:} Great question! We want to use the five tetrominoes\ldots
    \item \textbf{Reader (interrupting):} Tetrominoes?
    \item \textbf{Us:} Yes, that's the term we use for the pieces! ``Tetro-'' means ``four,'' and ``mino'' comes from ``domino,'' so this means ``like a domino but with four tiles''.
    \item \textbf{Reader:} Got it! Please continue.
    \item \textbf{Us:} Right! An area composed of a set of unit squares on the grid is considered ``enclosed'' if it is bounded by an edge-connected fence made of tetrominoes.
    \item \textbf{Reader:} By ``edge-connected,'' do you mean I can't connect the pieces by their corners?
    \item \textbf{Us:} More or less! While it’s possible for some pieces to touch at the corners, they still need to be connected through their edges. To clarify, the fence you build with the five tetrominoes must be rook-connected, meaning you should be able to move a chess rook from any square of the fence to any other square, as if the entire grid were a chessboard.
    \item \textbf{Reader (taking two tetrominoes):} I see! If I place two tetrominoes together like this (Figure~\ref{fig:corner_connect}, right), it is rook-connected, but not in this other arrangement (Figure~\ref{fig:corner_connect}, left), right?
    \item \textbf{Us:} Exactly! Now try to build a fence using four of the five pieces, say \polfont i, \polfont l, \polfont n, \polfont t.
    \item \textbf{Reader:} Hmm, like this (Figure~\ref{fig:tentative_tetra})?  Is that allowed?
    \item \textbf{Us:} Yes! And see, you’ve created two holes! It is allowed, but do you think it will give you the largest possible enclosed area?
    \item \textbf{Reader:} Probably not! Hmm, I’ll work on it a bit more and come back to you when I find a good solution using all five pieces.
\end{dialogue}

Puzzles, such as the one above, are a good way to engage people in thinking logically and using mathematics as a way to model them and solve them, but we wish to do so on a different level, finding ways to go beyond puzzles and really to ``gamify'' mathematical research to let it be experienced by the public. In its ideal form, it would also enable the participation of the public in the research itself, providing a sort of feedback loop between the researchers and the people experiencing the game, making it a collaborative endeavour.

The primary example of gamification of research presented in this paper, the Fence Challenge, was inspired by a generalization of the tetromino fence problem described above~\cite{LRMR23}. Additionally, we discuss a collaborative game, \textit{Le Carré du Diable}, developed with input from participants of the workshop \textit{Gamification of Research}, held during the conference \emph{\href{https://indico.cern.ch/event/1038038/}{Let's talk about outreach!}}.

During the same workshop, we also introduced two other fully developed instances of gamification of research: one based on the computational complexity of art-gallery problems on polycubes~\cite{AR21,LRMR22}, and another stemming from the study of higher-dimensional cubical sliding puzzles~\cite{BMRV23} (see \cref{ssec:Dominatrix,ssec:SlidingPuzzles}).

It is not the goal of this contribution to provide fully developed theoretical arguments regarding the benefits of gamification. The concept of gamification and its implications form a rich field of study that requires far more space to address comprehensively. Here, we aim, instead, to highlight certain design choices we made and share insights gained during our outreach activities. In the future, we hope to systematically analyse our setups and the behaviour of participants to better understand the impact of our approach.

Our outreach efforts aim to engage the widest audience possible, with particular emphasis on participants who might not typically enjoy mathematics. Beyond merely exposing them to mathematical concepts, we aim to immerse them in contemporary research and even invite them to contribute to it. Games serve as a powerful medium for achieving these goals by breaking down barriers and fostering curiosity. The ludic potential of games can transform attitudes toward mathematics~\cite{BHMW85}, and research shows that enjoyment of mathematics is one of the strongest predictors of student success~\cite{DDYF20}. Furthermore, games enable collective engagement, allowing participants and their surrounding communities to share the experience of mathematical discovery~\cite{MSJCJCSS17,QP22}.

The idea of presenting mathematics through games has a long and rich history. From the educational tools of the tutor in Rousseau’s \emph{Émile}~\cite{Rousseau1762} to Martin Gardner’s \emph{Mathematical Games} columns~\cite{GardnerMathGames}, games have been used to make mathematics accessible and engaging. Foundational works like Berlekamp, Conway, and Guy’s \emph{Winning ways for your mathematical plays}~\cite{BCG04} study games mathematically, inspiring further gamification efforts. Recent examples, such as Ravi Vakil’s ``bedtime story'' for understanding long exact sequences~\cite{Ra21}, Lee and Hua’s Homotopy Type Theory quests~\cite{HoTT}, and Buzzard and Pedramfar's \emph{Natural Number Game}~\cite{Buzz}, demonstrate how games can bridge the gap between recreational engagement and mathematical research. These examples illustrate the transformative potential of games, and we hope that the experiences shared here will inspire others to explore the interplay between games and research-level mathematics.

We now provide an overview of the sections in this contribution. Section~\ref{sec:gamification} briefly presents our first two examples of gamification along with some of our design choices. Section~\ref{sec:fence_challenge} focuses on the Fence Challenge, the gamification exercise proposed to participants of \emph{Let's talk about outreach!}, and the collaborative game developed during this workshop. In Section~\ref{sec:Conclu}, we synthesise our observations and outline the next steps in our research. Lastly, readers interested in the puzzle introduced earlier in this section can find its complete solution in \cref{app:Proof_Tetra}, as well as supplementary figures for the game in  \cref{app:Ex}.

\section{Two examples of gamification}\label{sec:gamification}

\subsection{DominatriX}\label{ssec:Dominatrix}

The activity ``DominatriX'' gamifies research on domination problems on polycubes, as explored in~\cite{AR21,LRMR22}. These papers show that the minimal domination problem of queens and rooks on polyominoes and the independent set problem of queens and rooks on polycubes of dimension $d \geq 3$ are NP-complete. Conversely, the independent set problem for rooks on polyominoes is solvable in polynomial time (P). The status of the independent set problem for queens on polyominoes remains unresolved but is conjectured to be NP-complete.

The activity aims to let participants experience the intrinsic difficulties of these domination problems through a video game and guided activities. Additionally, it can introduce advanced audiences to the distinction between P and NP complexity classes.

\paragraph{Setup} The activity is divided into three stations. The first station has chessboards with tokens representing queens and rooks. The centrepiece station features computers running the video game, with one connected to a projector. The final station includes large sheets of graph paper, pens, and tokens.  

\begin{figure}[t]
    \centering
\includegraphics[width=0.65\linewidth]{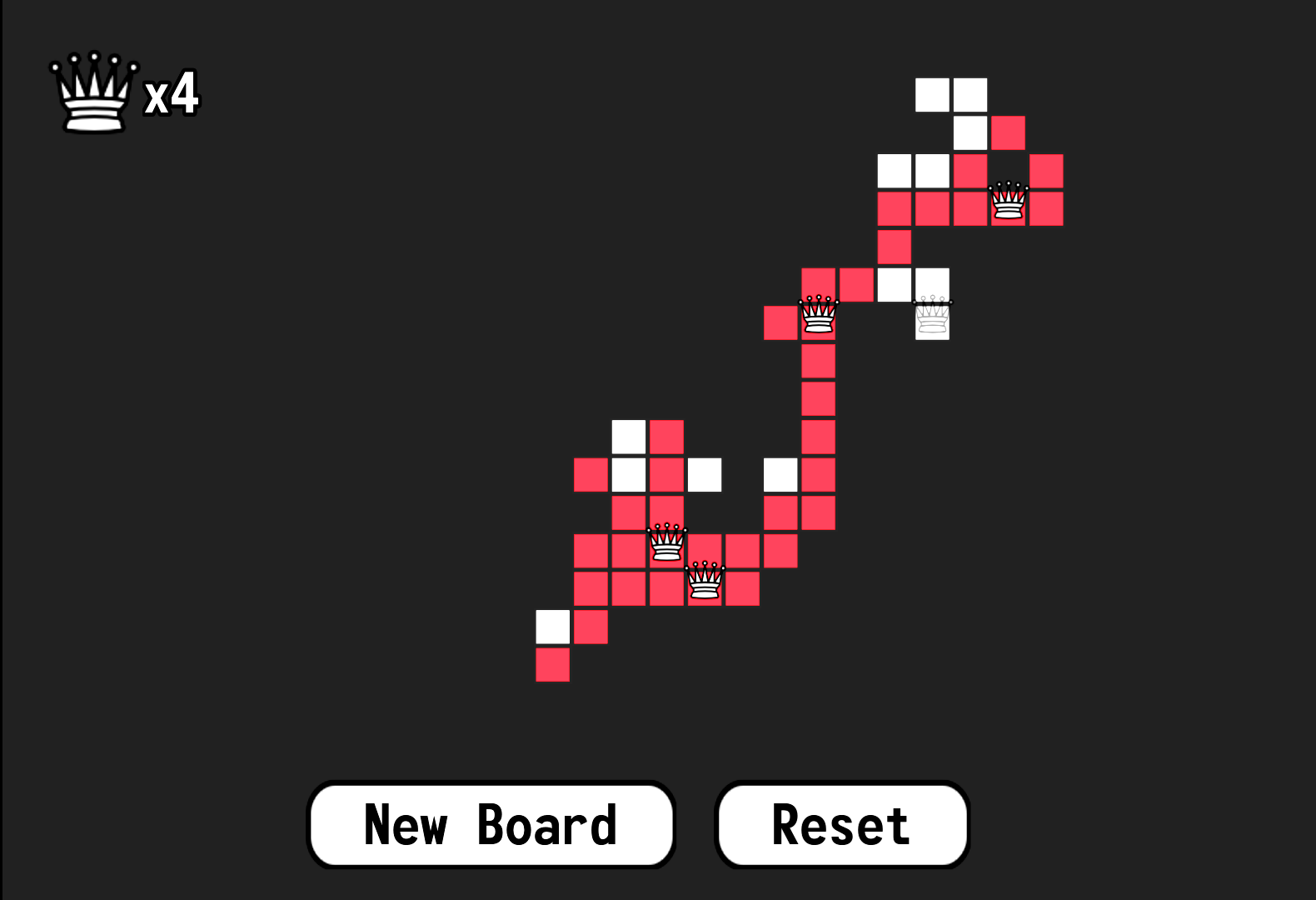}
    \caption{Screenshot from a level of the video game \emph{Queens and Rooks Domination} accessible at \href{https://www.erikaroldan.net/queensrooksdomination}{https://www.erikaroldan.net/queensrooksdomination}. (The solution is 7 queens.)}
    \label{fig:video_game_QR}
\end{figure}

\paragraph{Content} The activity consists of a video game, animations, and a hands-on construction station. The video game challenges players to solve the minimal queen or rook domination problem on randomly generated 50-tile polyominoes. Players place pieces on the board, lighting up tiles threatened by each piece, with the goal of illuminating all tiles using as few pieces as possible. Upon completing a challenge, players receive feedback on the minimum number of pieces required and can retry, switch pieces, or start a new challenge. See \cref{fig:video_game_QR} for an example. The video game is accessible at~\href{https://www.erikaroldan.net/queensrooksdomination}{https://www.erikaroldan.net/queensrooksdomination}.

The animation introduces the idea of domination using chessboards, starting with the classical $n$-queen problem. Participants are guided to dominate chessboards of increasing size (from $1 \times 1$ to $8 \times 8$) with as few queens as possible or to place as many queens or rooks as possible without any two pieces threatening each other.

At the last station, participants can draw their own polyominoes and use them as the basis for solving the challenges. For example, two teams can exchange polyominoes and attempt to solve each other's puzzles. Animators can use the solver from~\cite{LRMR22} to find optimal solutions for the polyominoes.

\paragraph{Goals for advanced audiences} This activity introduces participants to computational complexity concepts, particularly the distinction between ``easy'' and ``hard'' problems. For instance, the independent set problem for rooks on polyominoes is ``easy'' (P), while the minimum domination problem is ``hard'' (NP-complete). An ``easy'' problem is solvable by an algorithm whose runtime grows polynomially with the input size, while a ``hard'' problem lacks such an efficient solution, assuming $\mathrm{P} \neq \mathrm{NP}$. Participants can experience this distinction firsthand: the independent set problem has a straightforward winning strategy, while the minimal domination problem requires heuristics that may not work for all instances.

\subsection{Cubical Sliding Puzzles}\label{ssec:SlidingPuzzles}

Sliding puzzles were popularised in the late 19th century by the 15-Puzzle. In this puzzle, a $4\times 4$ square is tiled with 15 numbered tiles (from 1 to 15) and one empty space. The task is to return the tiles to an ordered configuration from a scrambled one. An interesting feature of this puzzle is that not all configurations are reachable from any position. This was famously highlighted by Samuel Loyd's bet: he challenged readers to find a sequence of moves that would transform the starting ordered configuration into one where the ``14'' and ``15'' tiles were swapped. Loyd's money was never at risk, however, as it was proven impossible by Johnson and Story over 10 years earlier~\cite{JS1879}!

This outreach experience focuses on higher-dimensional sliding puzzles, as studied in~\cite{BMRV23}, where the solvability regimes of many higher-dimensional generalizations of sliding puzzles are investigated.

\begin{figure}[t]
    \centering
    \includegraphics[width=0.65\linewidth]{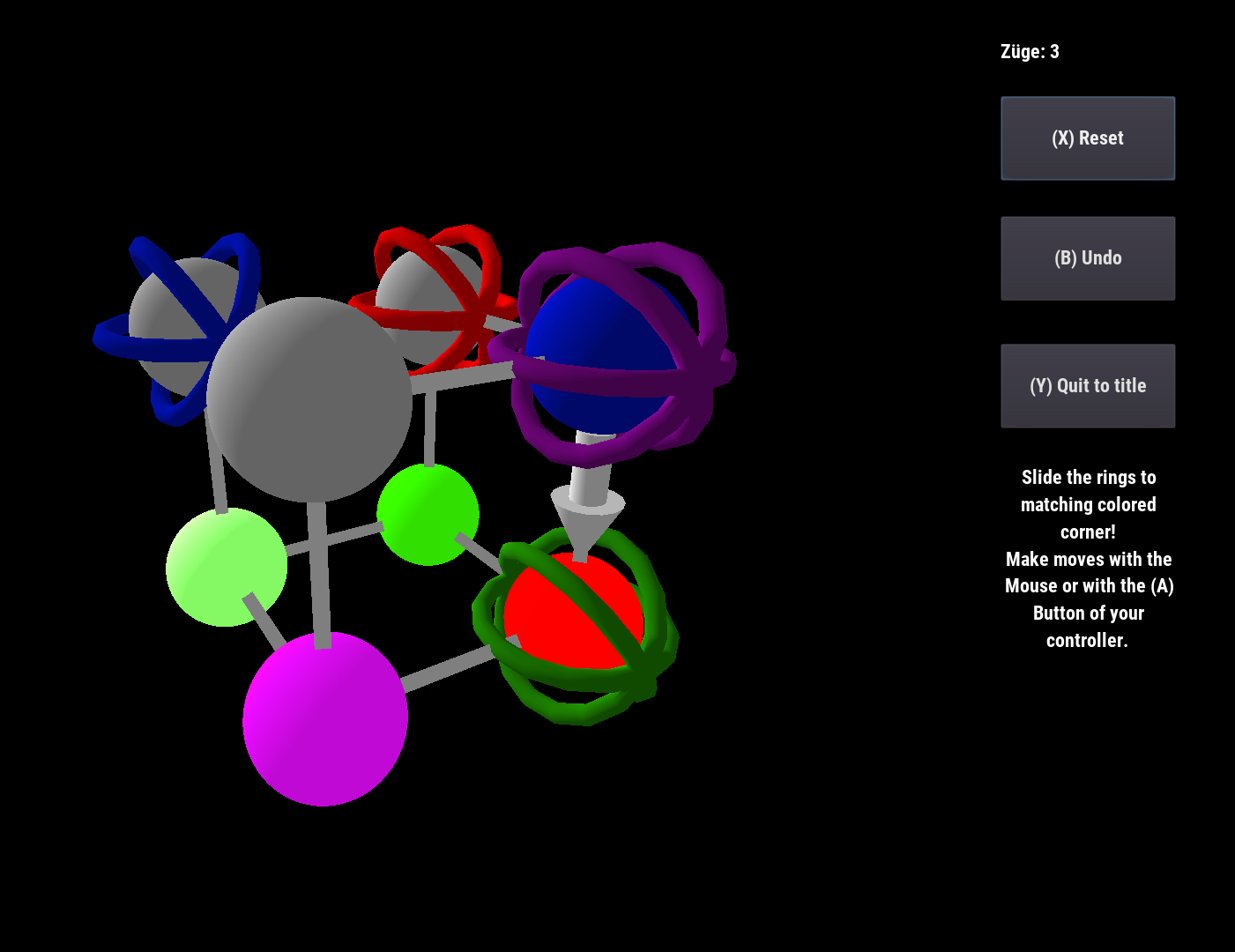}
    \caption{Screenshot of a cube-level from the video game \textit{Cubical Sliding Puzzle}, available at~\href{https://www.erikaroldan.net/cubical-sliding-puzzles}{https://www.erikaroldan.net/cubical-sliding-puzzles}. The highlighted nodes represent the possible moves for the selected green rings.}
    \label{fig:video_game_sliding}
\end{figure}

\paragraph{Setup} The activity setup includes a table with computers and game controllers connected to them (controllers are optional but often simplify movement). One computer is connected to a projector for larger group interactions. The table also features physical 2D sliding puzzles, such as the 15-puzzle or 3D-printed hexagonal sliding puzzles.

\paragraph{Content} The activity begins by allowing participants to explore the video game featuring cubical sliding puzzles (see \cref{fig:video_game_sliding}). Participants are encouraged to independently discover the game's mechanics and attempt to define the movement rules themselves. Notably, this is a challenging task in our experience. The video game is available at~\href{https://www.erikaroldan.net/cubical-sliding-puzzles}{https://www.erikaroldan.net/cubical-sliding-puzzles}.

In the second part of the activity, participants explore the relationship between these puzzles, the classic 15-Puzzle, and hexagonal sliding puzzles  \cite{karpman2022parity}. They play with the provided physical puzzles before progressing to four-dimensional cubical sliding puzzles in the video game. The challenge then becomes solving these puzzles in the fewest possible moves across the available levels.

\paragraph{Goals for advanced audiences} This game offers an opportunity to explore advanced mathematical concepts in topology, combinatorics, discrete configuration spaces, and the extension of graph-theoretic ideas to discrete topology through generalizations of the cell-cube structure in higher dimensions. Advanced participants first work to understand the movement rules in higher dimensions, which are governed by topological obstructions. They can then investigate the principles determining the solvability of the puzzles and work on finding solutions that minimize the number of moves.

\subsection{Design choices}\label{ssec:Gamification}

Here, we briefly outline the design choices made in crafting our activities, guided by our core values. Future research will be necessary to deepen this discussion once more data is collected and additional feedback from participants is received.

The first choice was to base the activities on the mathematics and computer science we actively research. While this approach poses challenges---such as the need to create new materials from scratch, the inability to answer all participant questions due to the open-ended nature of the research, the necessity to update activities as research evolves, and the difficulty in onboarding colleagues unfamiliar with the work---we believe the benefits far outweigh these drawbacks.

In fact, many of these ``drawbacks'' can be reframed as advantages. The novelty of the material ensures a unique experience for participants. Not knowing all the answers flattens the traditional hierarchy between presenter and participant, encouraging participants to propose ideas without fear of being ``wrong.'' Regular updates to reflect research progress offer recurring participants a dynamic view of mathematical discovery and enable us to refine the materials and activities. Introducing the activities to colleagues first serves as a valuable trial run, identifying potential improvements and ensuring clarity. Moreover, having facilitators who have only recently learned the activity is possible, and it could lower intimidation barriers, fostering a more inclusive environment while making them attuned to participants’ potential challenges.

The second choice was to design the activities to require minimal instruction, enabling them to be experienced with minimal guidance from us. This approach also aimed to reduce language barriers, allowing the games to be presented in various languages. To date, the activities have been conducted in different countries in America and Europe in Dutch, English, French, German, and Spanish, enabling us to engage participants from diverse linguistic and cultural backgrounds.

The third choice was to diversify the mediums through which participants interact with the activities. For example, in the queen domination activity, we include physical chessboards and pieces alongside computers or tablets. This ensures that participants who prefer not to, or cannot, use digital tools can still engage. Large screens projecting the games further allow bystanders to follow along and interact with active participants. Beyond accessibility, physical materials enhance comprehension by grounding abstract mathematical concepts in tangible objects. By physically engaging with the material, participants can embody mathematical ideas, gaining deeper insights through sensory feedback---a principle supported by research in embodied cognition~\cite{SABDD21,SBDD24}.

 We also aimed to establish a relationship of equality with participants. By providing opportunities to experiment independently---such as continuing to work on the activities through the freely accessible websites we have created---we encouraged participants to engage with the research process on their own terms rather than relying solely on us for guidance. This approach fosters collaboration and exploration, ensuring that participants feel empowered to continue actively contributing to the experience.
 
These two examples of gamification were designed to ensure that all our outreach encounters embody an inclusive culture and a strong commitment to welcoming the most diverse audience possible. Furthermore, they provided valuable insights and experience that informed the design of our third example of gamification of research, the Fence Challenge, which we present in the following section.

\section{The Fence Challenge and the game Le Carré du Diable} \label{sec:fence_challenge}
In this section, we will present succinctly the research questions considered in~\cite{LRMR23} that were the starting point of the gamification exercise proposed at the workshop and continue on the collaborative game the participants of the workshop designed. 

\subsection{Genesis of the problem: pentomino fence}

A set of pentomino fence puzzles was first presented by Feser in 1968~\cite{Fe68}. One of them, asking to find the greatest area enclosed by the 12 pentominoes, gained popularity when it appeared in the famous Gardner's \emph{Mathematical Games} column~\cite{Ga73}. In this column, Gardner presented his best attempt  and challenged the readers to find a solution matching the one found by Donald Knuth, who claimed to have a sketch of a proof that it was the maximal solution. (A solution is presented in \cref{fig:solmax}; the reader wanting to try by themselves might want to pause here and try on their own before we give the answer.)

The problem would often appear in the following editions, with readers claiming to have found proofs or giving some maximum values. There were many wrong attempts, some easily disprovable, since they tried to show the maximum to be smaller than the solution provided by Knuth. The first complete proof of the maximality we have found is the one of Takakazu Shimauchi~\cite{Ta78}. Sadly (for non-Japanese-speaking people) it was only available in Japanese. It is translated and presented with small adaptations in~\cite{LRMR23}.

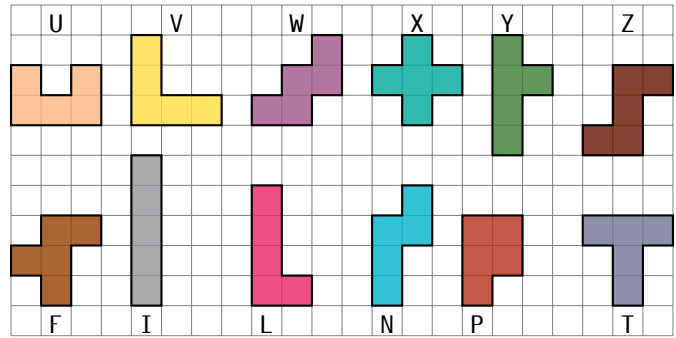
\begin{figure}[h]
    \centering
    \begin{tikzpicture}[scale=\scfa]
    \draw[help lines] (0,-1) grid (22,10);
        \pic[scale=\scfa] at (0,0) {pentominoF};
        \draw (1.5,0) node[below] {\polfont F};
        \pic[scale=\scfa] at (4,0) {pentominoI};
        \draw (4.5,0) node[below] {\polfont I};
        \pic[scale=\scfa] at (8,0) {pentominoL};
        \draw (8.5,0) node[below] {\polfont L};
        \pic[scale=\scfa] at (12,0) {pentominoN};
        \draw (12.5,0) node[below] {\polfont N};
        \pic[scale=\scfa] at (15,0) {pentominoP};
        \draw (15.5,0) node[below] {\polfont P};
        \pic[scale=\scfa] at (19,0) {pentominoT};
        \draw (20.5,0) node[below] {\polfont T};
        \pic[scale=\scfa] at (0,6) {pentominoU};
        \draw (1.5,10) node[below] {\polfont U};
        \pic[scale=\scfa] at (4,6) {pentominoV};
        \draw (5.5,10) node[below] {\polfont V};
        \pic[scale=\scfa] at (8,6) {pentominoW};
        \draw (9.5,10) node[below] {\polfont W};
        \pic[scale=\scfa] at (12,6) {pentominoX};
        \draw (13.5,10) node[below] {\polfont X};
        \pic[scale=\scfa] at (16,5) {pentominoY};
        \draw (16.5,10) node[below] {\polfont Y};
        \pic[scale=\scfa] at (19,5) {pentominoZ};
        \draw (20.5,10) node[below] {\polfont Z};
    \end{tikzpicture}    
    \caption{All the twelve pentominoes.}\label{fig:pentominoes}
\end{figure}

\subsection{New results on extremal fence problems}
We begin with the formal definition of the problem at hand. Recall that there are three tessellations of the plane: the triangular, the square and the hexagonal. We consider polyamonds, polyominoes, and polyhexes of a given size $n$ that we call the \emph{poly-elements} of the fence.

\begin{definition}
 Given a tessellation of the plane and poly-ements of size $n$, we call a fence a polyform constituted of all the poly-elements of size $n$ such that every element has two neighbours.
\end{definition}

This definition allows for many trivial fences. To know the board size required to place any optimal solution, or near-optimal solution, we can add  one further condition. Basically, it says it should enclose a single area with no spurious polyforms creating a different topology than a circle. It is a natural assumption for the polyomninoes enclosing the maximal area. 

\begin{lemma}\label{lem:BoardSize}
    If the pentomino fence is isotopic to the circle, then it can be placed on a $20\times 20$ board.
\end{lemma}

Another point of the research~\cite{LRMR23} is to give the first English proof of this problem (a translation and adaptation of the proof of Shimaushi as found in~\cite{Ta81}), and extend it to general settings of other tessellations.  

\begin{proposition}\label{prop:Pento_Farm}
The pentomino fence problem has 1440 solutions with optimal enclosed area of 128 tiles.
\end{proposition}

The proof of this proposition can be found in~\cite{LRMR23}, and we will not repeat it here. The idea of the proof is to bound the possible solutions by examining the largest rectangles that can contain the configurations while respecting a given perimeter. Then, the specific characteristics of the pieces are used to provide an upper bound for the maximal solution. It remains to establish a lower bound for the maximal solution by identifying candidate solutions. This step is handled through an algorithmic search. A similar result for the tetromino fence problem is provided in Appendix~\ref{app:Proof_Tetra}.

The final part of the research presented in~\cite{LRMR23} involved creating an efficient Integer Linear Programming (ILP) solver for this game. Due to the vast solution space and the rarity of optimal solutions, such a solver cannot be expected to find all optimal solutions for fences made of larger poly-elements. However, for smaller sizes, it is capable of finding all optimal solutions. For instance, it can determine optimal solutions for subsets of pentominoes, which presents an interesting challenge.

The solver is versatile, working for any lattice in any dimension. Additionally, the research in~\cite{LRMR23} explores related questions for other tessellations and polyforms.

\subsection{The game -- Le Carré du Diable}\label{sec:5Game}

In this section, we present the game co-created with the participants of the workshop, along with some variants. The objective of the game is to collectively construct the fence enclosing the biggest possible area with the pentomino using a limited number of moves, starting from a given initial position. Players take turns sequentially, and the game concludes either when the allotted turns are exhausted or when all players choose to pass.

\subsection{Material and setup}
This game is played on a $20\times 20$ board using the 12 pentominoes (as shown in \cref{fig:pentominoes}). The board size is derived directly from \cref{lem:BoardSize}. While the game only requires a board and the pentominoes, we found it useful to incorporate a mixed reality component using a web-based app and additional materials, which are described in \cref{ssec:MixedReality}.

It is worth noting that coincidentally, a $20\times 20$ board is the same size as the board used in the game Blokus~\cite{Blokus}. This provides a cheap and practical way to play: simply use a Blokus board and a set of 12 pentominoes. From experience, we recommend mixing the colours of the pentominoes to make it easier to distinguish them visually.

\subsection{Rules}

We now describe how to play the collaborative game created at the workshop. A set of detailed examples are presented in Appendix~\ref{app:Ex}.

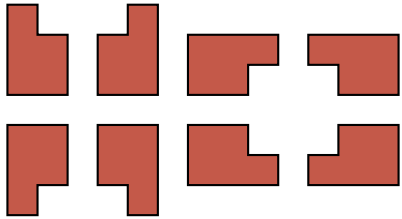
\begin{figure}[h]
\centering
 \begin{tikzpicture}[scale = \scfa]
    \pic[scale=\scfa] at (0,0) {pentominoP};
    \pic[scale=\scfa] at (5,0) {pentominoP={-1}};
    \pic[scale=\scfa, rotate = 90] at (9,1) {pentominoP};
    \pic[scale=\scfa, rotate = 270] at (10,1) {pentominoP={-1}};
    \pic[scale=\scfa, rotate = 180] at (5,7) {pentominoP};
    \pic[scale=\scfa, rotate = 90] at (9,6) {pentominoP={-1}};
    \pic[scale=\scfa, rotate = 270] at (10,6) {pentominoP};
    \pic[scale=\scfa, rotate = 180] at (0,7) {pentominoP={-1}};
    \end{tikzpicture}
        \caption[All the ways to place the piece ``P''.]{All the ways to place a piece, here the piece ``\begin{tikzpicture}[scale=.1,baseline={(current bounding box.center)}]
            \pic[scale=.1] at (0,0) {pentominoP};
    \end{tikzpicture}''.}
    \label{fig:refl_rot_F}
    \end{figure}

\paragraph{General rules}
\begin{enumerate}
    \item One person plays one move at a time.
    \item The pieces must always be placed on the grid.
    \item All rotations and reflections of the pieces are allowed (see Figure~\ref{fig:refl_rot_F}).
    \item Each piece must be edge-connected to at least two neighbours. It is important to note that being connected on a corner does \textbf{not} count as being connected (see Figure~\ref{fig:ex_coup} for an example and a non-example).
    \item Each player makes $\lceil 24\div \text{(number of participants)}\rceil $ moves\footnote{We recommend simply providing a table like Table~\ref{table:NumberMove}.}.
\end{enumerate}

\vspace{.5cm}

\begin{minipage}{0.45\textwidth}
\paragraph{Course of a game}
\begin{enumerate}
    \item The players choose a starting position.
    \item The person who made a multiplication most recently starts. Then, the turns continue clockwise.
    \item A participant can always pass their turn, but this counts in the number of moves at their disposition.
    \item The game ends when everyone passes or when everyone has made all their moves.
    \item Count the number of tiles inside the fence, more is better. 
\end{enumerate}
\end{minipage}
\hspace{0.05\textwidth}
\begin{minipage}{0.45\textwidth}
\paragraph{Course of a move}
\begin{enumerate}
    \item Choose to play or to pass. If the turn is passed, then the move goes to the next player in clockwise order.
    \item Choose a piece and move it, it is possible to flip and rotate it.
    \item Make sure the board is valid: each piece must have 2 neighbours. If the board is invalid, the piece is put back to its original place and the move starts again.
    \item When the move is played, the turn passes to the next person in clockwise order.
\end{enumerate}
\end{minipage}

\paragraph{Notes on the game}
\begin{enumerate}
    \item You have to play only one piece; you cannot exchange two pieces.
    \item Attention, a move is not forced to improve the area. Even, sometimes, a move can reduce the area and prepare the way for a better position.
    \item A score of 100 or more is good, 120 or more is excellent, 125 or more is exceptional, and 128 is the maximum. 
\end{enumerate}

\begin{table}[hp]
    \centering
\begin{tabular}{l|rrrrrrrrrrrr}
Number of people & 1  & 2  & 3 & 4 & 5 & 6 & 7 & 8 & 9 & 10 & 11 & 12 \\ \hline
Number of moves    & 24 & 12 & 8 & 6 & 5 & 4 & 3 & 3 & 3 & 2  & 2  & 2 
\end{tabular}
    \caption{Number of moves by person.}
    \label{table:NumberMove}
\end{table}
\subsection{Variants}
Through the interactions with the participants some variants of the game were suggested to us. We list them here with a small explanation in the case someone would like to pick them up.

\paragraph{Smaller sets of pieces}
This variant, good for quicker games or for younger participants, is to take a subset of the 12 pentominoes. 

\paragraph{Misère}
In this version, the players try to reduce as much as possible the size of the hole(s) in a given set of moves. The goal is to reach a closed polyshape.

\paragraph{Hexominoes}
The same game is given, but this time with the hexominoes. This is quite big, since there are 35 hexominoes.

\paragraph{Pentahexes}
The same game is given on a hexagonal grid with the 22 pentahexes. This game is also quite hard, and the maximum is not known. We have constructed one version with laser-cut pieces placed on a role-playing game mat. Here, counting the number of hexagons enclosed by the fence is quite cumbersome. 

\subsection{Mixed reality}\label{ssec:MixedReality}
Keeping on the idea of accessibility, we have created a program that computes the area of a given configuration of pentominoes. This allows for participants to immediately get feedback on their moves without having to count the whole area enclosed. This plays into the relation between the object and the impacts of a move by giving feedback on the action of the participants: moving a piece has consequences that can be immediately seen, and as such, it helps to quickly build up the intuition. It is available at this webpage \href{https://www.erikaroldan.net/fencechallenge}{https://www.erikaroldan.net/fencechallenge}.

\subsection{Setup of one outreach experience}

We now describe an outreach experience inspired by this game, as presented to participants during the ``Lange Nacht der Wissenschaften'' (Long Night of the Sciences) on June 23rd, 2023, in Leipzig, organised and presented by the first and third authors on behalf of the Max Planck Institute for Mathematics in the Sciences.

\paragraph{Material} 
Our stand was part of the activities hosted at the Max Planck Institute for Mathematics in the Sciences, alongside various other mathematical exhibits. The setup consisted of a table with the game grid, a camera, lighting, and a screen positioned above the stand to display the activity.

\paragraph{Places and events} 
The stand was active from 18:00 to 22:00 in one of the two rooms designated for the evening’s activities. Participants could explore many other games alongside ours. Each stand was staffed by one accompanying person. The room also featured other gamified research activities, such as a Mondrian-type puzzle~\cite{GCLMR23}, a labyrinth, and the chess-domination puzzle~\cite{LRMR23}, all of which were developed and implemented by students and researchers from the third author’s research groups. Additional activities were available outdoors. 

Two scheduled animations presented the mathematical ideas behind the games. Participation in the event was free, and attendees were encouraged to explore as many stands as they wished. Playing three games earned participants a token that could be exchanged for a small mathematics-themed gift.

\paragraph{Course of the activity} 
The instructions provided were intentionally minimal, partly due to the constraints of the volunteer team---most of whom had only a basic command of German---and partly by design, as we wanted participants to discover the rules independently.

The instructions were typically divided into three phases. In the first phase, participants were given an initial position that could be easily improved (see Figure~\ref{fig:game_example}\,(a)). Participants could then either interact with the game, using computer feedback to improve their score, or have the rules explained: how to place pieces, the goal of increasing the area, etc. 

In the second phase, small groups of participants played a shortened version of the game described in Section~\ref{sec:5Game}, usually lasting 10--12 moves. In the final phase, participants were encouraged to freely experiment with the game, attempting to improve the size of the solution and beat the highest score\footnote{Interestingly, a bug in the software caused it to crash for solutions of 126 tiles or more. This ``feature'' delighted the few participants who managed to reach and exceed this threshold!}.

\section{Conclusion}\label{sec:Conclu}

We believe the three examples presented here, along with the workshop we designed and implemented, demonstrate the potential of gamifying research in mathematics and computer science. During our presentations, we observed that these activities effectively engaged participants in mathematical thinking, sparked curiosity, and were widely appreciated by attendees.

The gamification of research serves not only as a tool to popularise new findings in mathematics and computer science, but also as a means to inspire curiosity and explore new research directions. At its core, we define the gamification of research as \emph{a methodology to develop, curate, and implement outreach and educational projects that allow participants to explore the concepts underlying mathematical research in a playful and engaging manner}. The target audience can vary, ranging from the general public, as in the chess domination or pentomino fence games, to specialised groups, such as mathematicians learning specific tools through the HoTT and Natural Number games.

Our work adheres to the following guiding principles:
\begin{enumerate}
    \item \textbf{Accessibility.} Activities are designed to be inclusive and enjoyable for as many people as possible.
    \item \textbf{Originality.} Each activity addresses a subject currently under investigation, ensuring novelty.
    \item \textbf{Visualisation of mathematics and choice of media.} Diverse formats and tools are employed to effectively engage participants.
    \item \textbf{Horizontality.} We strive to break traditional vertical hierarchies of knowledge transmission, encouraging collaborative and participatory learning.
\end{enumerate}

A key aspect of gamification that we wish to emphasize is the feedback loop it can create with participants. For instance, the ILP program presented in~\cite{LRMR22} was developed in response to design challenges encountered during the creation of the chess domination game, which itself was inspired by a research question posed in~\cite{AR21}.

\subsection{Future research}

While some of us have been involved in the gamification of mathematical research for many years, our theoretical exploration of the concept is still relatively new. Based on our experience, we believe gamification to be a deeply enriching endeavour for both participants and our academic practices.

Looking ahead, we aim to analyse the theoretical foundations of gamification by engaging with recent literature in game studies. This work will be further supported by a planned quantitative investigation of research experiences and games developed within this framework, leveraging educational technology and learning analytics research methodologies.

Promising directions for future exploration include:
\begin{enumerate}
    \item Assessing the impact of games and activities on participants' perceptions of mathematics over the short, medium, and long term.
    \item Investigating how participant feedback contributes to and shapes subsequent mathematical research.
    \item Evaluating the effectiveness of games in engaging diverse audiences, including students at different academic levels, STEAM researchers, and the general public.
    \item Developing systematic strategies to reach and create participatory STEAM opportunities for marginalised audiences by consciously designing inclusive materials and activities. 
\end{enumerate}

\section*{Acknowledgements}
We acknowledge the contributions of Peter Voran and Mia Mü\ss{}ig with whom the first- and third-named co-authors animated and co-created the workshop of \textit{Gamification of Research}. We wish to thank the participants of the conference \emph{Let’s talk about outreach!} that took place at the SwissMAP Research Station, Les Diablerets, Switzerland for their help creating the fence game. In particular,  we thank Hacène Belbachir, Pierre-Alain Cherix, Shaula Fiorelli, Timothée Jaccard, Samuel Lelièvre, Ayliean McDonald, Cesco Reale and Alain Valette, who devised the first rules of the cooperative game \emph{Le Carré du Diable}.
We also wish to thank Johannes Häfner for his help creating the mixed reality program. ALR also wishes to thank Ken Lee and Joseph Hua for accepting to discuss their initiative with him. Final thanks go to Tamara Sprinkle for her professional and thorough translation of the article~\cite{Ta78}.

ALR also acknowledges the hospitality of the Max Planck Institute for Mathematics in the Sciences during his visits and ScaDS.AI Leipzig where he was employed during part of this research.

\printbibliography

\appendix

\section{Proof of the tetromino fence problem}\label{app:Proof_Tetra}

In this appendix, we solve the introduction problem. This presents the idea of the proof of Proposition~\ref{prop:Pento_Farm}.
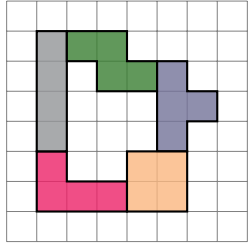
\begin{figure}[hp]
    \centering
    \begin{tikzpicture}[scale=\scfa]
    \draw[help lines] (0,0) grid (8,8);
        \pic[scale = \scfa,rotate=270] at (1,1){tetrominoL={-1}};
        \pic[scale=\scfa,rotate=0] at (1,3) {tetrominoI};
        \pic[scale = \scfa,rotate =0] at (5,5) {tetrominoN={-1}};
        \pic[scale = \scfa,rotate=90] at (7,3){tetrominoT};
        \pic[scale=\scfa,rotate=0] at (4,1) {tetrominoO};
    \end{tikzpicture}
    \caption{One of the optimal 9-tile solutions to the tetromino fence problem.}
    \label{fig:sol_tetra}
\end{figure}

\begin{proposition}
The tetromino fence problem has 21 solutions of maximal area 9.
\end{proposition}
\begin{proof}
    Let us begin by computing the biggest perimeter possible using the  five tetrominoes of~\eqref{eq:tetrapieces}. For this, we will analyse the \emph{length} of the pieces in the sense of by how much can they increase the perimeter. For clarity purpose, we will introduce some vocabulary. Given a tetromino, we first identify the \emph{direction} of the tetromino, that is, where is it the longest. Then, we take an initial tile at an extremity of the tetromino, and  we compute two quantities: first the \emph{progression}, that is, how many more squares do the tetromino have in the direction of progress; second, the \emph{protrusion}, that is, how many squares do the tetromino have in the other direction of progress. The length is then the progression added to the protrusion plus the initial tile. 
    Some pieces might be employed differently with the same maximal length. Let us go through each five pieces now.
\begin{itemize}
    \item The ``\polfont{i}'' tetromino is the simplest, it progresses straight for 3 more tiles. Progression: 3, protrusion: 0, length: 4.
    \[
    \begin{tikzpicture}[scale = \scfa]
        \pic[scale=\scfa,rotate=90] at (0,0) {tetrominoI};
         \foreach \x/\y in {-4/0,-3/0,-2/0,-1/0}
        {
        \fill (\x+.5,\y+.5) circle(2pt);
        }
        \draw[->] (-3.5,.5) -- (-.5,.5);
\end{tikzpicture}
    \] 
    \item The ``\polfont{l}'' tetromino progresses straight for 2 more tiles and then bends for 1. Progression: 2, protrusion: 1, length: 4.
    \[
    \begin{tikzpicture}[scale = \scfa]
        \pic[scale=\scfa,rotate=90] at (0,0) {tetrominoL};
        \foreach \x/\y in {-3/0,-2/0,-1/0,-1/1}
        {
        \fill (\x+.5,\y+.5) circle(2pt);
        }
        \draw[->] (-2.5,0.5) -- (-.5,.5);
        \draw[->] (-.5,.5) -- (-.5,1.5);
\end{tikzpicture}
    \] 
    \item The ``\polfont{n}'' tetromino progresses straight for 1 more tile and then bends for 1 and then progresses of 1. Progression: 2, protrusion: 1, length: 4.
    \[
    \begin{tikzpicture}[scale = \scfa]
        \pic[scale=\scfa] at (0,0) {tetrominoN};
        \foreach \x/\y in {0/0,1/0,1/1,2/1}
        {
        \fill (\x+.5,\y+.5) circle(2pt);
        }
        \draw[->] (.5,0.5) -- (1.5,.5);
        \draw[->] (1.5,.5) -- (1.5,1.5);
        \draw[->] (1.5,1.5) -- (2.5,1.5);
\end{tikzpicture}
    \] 
    \item The ``\polfont{t}'' tetromino can be used in two ways. One way is to make it progress straight for 2 more tiles. Alternatively, it could also progress straight for one and then bends for 1. In both alternatives, it has length 3. A) Progression: 2, protrusion: 0, length: 3; B) Progression: 1, protrusion: 1, length: 3.
    \[
    \begin{tikzpicture}[scale = \scfa]   
        \pic[scale=\scfa,rotate=180] at (0,0) {tetrominoT};
        \draw (-1,-2) node[below] {A)};
        \foreach \x/\y in {-3/-2,-2/-2,-1/-2,-2/-1}
        {
        \fill (\x+.5,\y+.5) circle(2pt);
        }
        \draw[->] (-2.5,-1.5) -- (-0.5,-1.5);
        \pic[scale=\scfa,rotate=180] at (10,0) {tetrominoT};
        \draw (9,-2) node[below] {B)};
        \foreach \x/\y in {-3/0,-2/0,-1/0,-2/1}
        {
        \fill (\x+.5+10,\y+.5-2) circle(2pt);
        }
        \draw[->] (-3.5+11,-1.5) -- (-2.5+11,-1.5);
        \draw[->] (-2.5+11,-1.5) -- (-2.5+11,-.5);
\end{tikzpicture}
    \] 
    \item Finally, the ``\polfont{o}'' tetromino progresses of 1 and then bends for 1. Progression: 1, protrusion: 1, length: 3.
    \[
    \begin{tikzpicture}[scale = \scfa]
        \pic[scale=\scfa] at (0,0) {tetrominoO};
        \foreach \x/\y in {-0/0,1/0,0/1,1/1}
        {
        \fill (\x+.5,\y+.5) circle(2pt);
        }
        \draw[->] (0.5,0.5) -- (1.5,.5);
        \draw[->] (1.5,.5) -- (1.5,1.5);
\end{tikzpicture}
    \] 
\end{itemize}
Adding the lengths together, we reach a possible perimeter of 18. This can be achieved by a rectangle of $5\times 4$, and such a rectangle has an inside area of 12. However, due to the geometry of the pieces \polfont{o} and \polfont n, there will always be 3 tiles from the perimeter inside the inner area: 1 from \polfont o and 2 from \polfont n when they are used as optimal length. Hence, 9 is the maximal inner area a shape can achieve when the tetrominoes are circumscribed by a $5\times 4$ rectangle. One such solution is given in Figure~\ref{fig:sol_tetro}. 
\begin{figure}[h]
    \centering
    \begin{tikzpicture}[scale=\scfa]
    \draw[help lines] (0,0) grid (8,8);
        \pic[scale = \scfa,rotate=270] at (1,1){tetrominoL={-1}};
        \pic[scale=\scfa,rotate=0] at (1,3) {tetrominoI};
        \pic[scale = \scfa,rotate =0] at (5,5) {tetrominoN={-1}};
        \pic[scale = \scfa,rotate=90] at (7,3){tetrominoT};
        \pic[scale=\scfa,rotate=0] at (4,1) {tetrominoO};
        \draw[thick,dotted] (1.5,1.5)  rectangle (5.5,6.5);
        \draw[fill=red,fill opacity=.75] (4,2) rectangle (5,3) ;
        \draw[fill=red,fill opacity=.75] (3,6) rectangle (5,5) ;
        \foreach \x/\y in {2/2,3/2,2/3,3/3,4/3,2/4,3/4,4/4,2/5}
        {
        \fill (\x+.5,\y+.5) circle(2pt);
        }
    \end{tikzpicture}
    \caption{One of the optimal 9-tile solutions to the tetromino fence problem. The circumscribing rectangle is drawn by a dotted line, and the defects inside are shown in shaded red. The inner tiles have a dot on them}\label{fig:sol_tetro}
\end{figure}
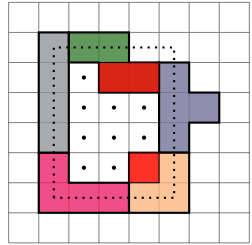

The other possible rectangles circumscribing the solution of perimeter less than 18 are $4\times 4$ (perimeter: 16) and $5\times 3$ (perimeter: 16), $3\times 4$ (perimeter: 12). The inside areas are 9, 8 and 6, respectively. Hence, the maximum possible solution is indeed 9.
    
    The 21 solutions were then found by brute force, looking at all possible configurations in the $5\times 4$ and $4\times 4$ rectangles\footnote{Computational results presented in~\cite{LRMR23}.}.
\end{proof}

\section{Examples}\label{app:Ex}

This appendix includes extra figures giving more explantations for the game.

\begin{figure}[hp]
    \centering
    \begin{tikzcd}
     \begin{tikzpicture}[scale = \scfa,baseline={(current bounding box.center)}]
\pic[scale=\scfa,rotate=270] at (0,3) {pentominoV={1}};
\pic[scale=\scfa] at (-1,3) {pentominoX={1}};
\pic[scale=\scfa,rotate= 90] at (6,5) {pentominoL ={-1}};
\pic[scale=\scfa,rotate= 270] at (2,1) {pentominoN ={-1}};
\end{tikzpicture}
    \arrow[rr,OliveGreen] \arrow[d,OliveGreen] \arrow[drr, "/" marking,red]&&
    \begin{tikzpicture}[scale = \scfa,baseline={(current bounding box.center)}]
\pic[scale=\scfa,rotate=0] at (-1,1) {pentominoV={1}};
\pic[scale=\scfa] at (-1,3) {pentominoX={1}};
\pic[scale=\scfa,rotate= 90] at (6,5) {pentominoL ={-1}};
\pic[scale=\scfa,rotate= 270] at (2,1) {pentominoN ={-1}};
\end{tikzpicture}\\
    \begin{tikzpicture}[scale = \scfa,baseline={(current bounding box.center)}]
\pic[scale=\scfa] at (0,0) {pentominoV={1}};
\pic[scale=\scfa] at (-1,3) {pentominoX={1}};
\pic[scale=\scfa,rotate= 90] at (6,5) {pentominoL ={-1}};
\pic[scale=\scfa,rotate= 270] at (2,1) {pentominoN ={-1}};
\end{tikzpicture}&&
    \begin{tikzpicture}[scale = \scfa,baseline={(current bounding box.center)}]
\pic[scale=\scfa,rotate=270] at (-1,3) {pentominoV={1}};
\pic[scale=\scfa] at (-1,3) {pentominoX={1}};
\pic[scale=\scfa,rotate= 90] at (6,5) {pentominoL ={-1}};
\pic[scale=\scfa,rotate= 270] at (2,1) {pentominoN ={-1}};
\draw[thick,dotted,red] (1-.1,3+.1) -- (3+.1,3+.1) -- (3+.1,1-.1) -- (1-.1,1-.1) -- cycle;
\end{tikzpicture}
    \end{tikzcd}
    \caption[We wish to move the yellow ``V'' piece. Examples of valid moves on the bottom left and top right. Example of an invalid move at the bottom right (the ``V'' and ``Z'' pieces are connected only by an edge).]
    {We wish to move the yellow ``\begin{tikzpicture}[baseline={(current bounding box.center)}]
\pic[scale=0.1,rotate=270] at (0,0) {pentominoV};
\end{tikzpicture}'' piece. Examples of valid moves on the bottom left and top right. Example of an invalid move at the bottom right (the fence is not closed: the yellow ``\begin{tikzpicture}[baseline={(current bounding box.center)}]
\pic[scale=0.1,rotate=270] at (0,0) {pentominoV};
\end{tikzpicture}'' and blue ``\begin{tikzpicture}[baseline={(current bounding box.center)}]
\pic[scale=0.1,rotate=270] at (0,0) {pentominoN={-1}};
\end{tikzpicture}'' pieces are only corner-connected, not edge-connected.}
    \label{fig:ex_coup}
\end{figure}
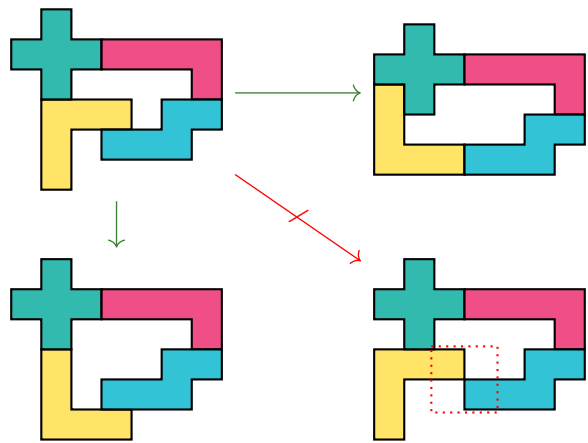

\begin{figure}[hp]
\begin{subfigure}{.45\textwidth}
    \centering
     \begin{tikzpicture}[scale=\scfa]
      \draw[help lines] (4,-1) grid (18,12);
         \pic[scale=\scfa,rotate=270] at (6,4) {pentominoI={-1}};
         \pic[scale=\scfa] at (8,4) {pentominoL};
         \pic[scale=\scfa,rotate=90] at (8,1) {pentominoT};
         \pic[scale=\scfa] at (8,7) {pentominoF};
        \pic[scale=\scfa] at (7,0) {pentominoX};
        \pic[scale=\scfa,rotate=180] at (14,10) {pentominoV};
        \pic[scale=\scfa] at (13,4) {pentominoN};
        \pic[scale=\scfa,rotate=180] at (12,3) {pentominoP};
        \pic[scale=\scfa,rotate=270] at (12,3) {pentominoU};
        \pic[scale=\scfa] at (14,0) {pentominoW};
        \pic[scale=\scfa] at (16,2) {pentominoY={-1}};
        \pic[scale=\scfa] at (5,4) {pentominoZ};
     \end{tikzpicture}
    \caption{Starting point: a fence enclosing 29 tiles. With three people playing, each will play 8 turns.}
    \end{subfigure}
    \hfill
    \begin{subfigure}{.45\textwidth}
    \centering
     \begin{tikzpicture}[scale=\scfa]
      \draw[help lines] (4,-1) grid (18,12);
         \pic[scale=\scfa,rotate=270] at (6,4) {pentominoI={-1}};
         \pic[scale=\scfa] at (8,4) {pentominoL};
         \pic[scale=\scfa,rotate=90] at (8,1) {pentominoT};
         \pic[scale=\scfa] at (8,7) {pentominoF};
        \pic[scale=\scfa] at (7,0) {pentominoX};
        \pic[scale=\scfa,rotate=180] at (14,10) {pentominoV};
        \pic[scale=\scfa] at (16,6) {pentominoN={-1}};
        \pic[scale=\scfa,rotate=180] at (12,3) {pentominoP};
        \pic[scale=\scfa,rotate=270] at (12,3) {pentominoU};
        \pic[scale=\scfa] at (14,0) {pentominoW};
        \pic[scale=\scfa] at (16,2) {pentominoY={-1}};
        \pic[scale=\scfa] at (5,4) {pentominoZ};
     \end{tikzpicture}
    \caption[First move: Andréa moves `` N ''. ]{First move: Andréa moves ``\begin{tikzpicture}[baseline={(current bounding box.center)}]
\pic[scale=0.1] at (0,0) {pentominoN};
\end{tikzpicture}'' and flips it. The fence now encloses 34 tiles.}
    \end{subfigure}
        \\[1cm]
    \begin{subfigure}{.45\textwidth}
    \centering
     \begin{tikzpicture}[scale=\scfa]
      \draw[help lines] (4,-1) grid (18,12);
         \pic[scale=\scfa,rotate=0] at (7,7) {pentominoI={-1}};
         \pic[scale=\scfa] at (8,4) {pentominoL};
         \pic[scale=\scfa,rotate=90] at (8,1) {pentominoT};
         \pic[scale=\scfa] at (8,7) {pentominoF};
        \pic[scale=\scfa] at (7,0) {pentominoX};
        \pic[scale=\scfa,rotate=180] at (14,10) {pentominoV};
        \pic[scale=\scfa] at (16,6) {pentominoN={-1}};
        \pic[scale=\scfa,rotate=180] at (12,3) {pentominoP};
        \pic[scale=\scfa,rotate=270] at (12,3) {pentominoU};
        \pic[scale=\scfa] at (14,0) {pentominoW};
        \pic[scale=\scfa] at (16,2) {pentominoY={-1}};
        \pic[scale=\scfa] at (5,4) {pentominoZ};
     \end{tikzpicture}
    \caption[Second move: Billie moves the ``I'' tile. ]{Second move: Billie moves the ``\begin{tikzpicture}[baseline={(current bounding box.center)}]
\pic[scale=0.1,rotate=90] at (0,0.5) {pentominoI};
\end{tikzpicture}'' tile and turns it 90\textdegree. The fence encloses 40 tiles.}
    \end{subfigure}
    \hfill
        \begin{subfigure}{.45\textwidth}
    \centering
     \begin{tikzpicture}[scale=\scfa]
      \draw[help lines] (4,-1) grid (18,12);
         \pic[scale=\scfa,rotate=0] at (7,7) {pentominoI={-1}};
         \pic[scale=\scfa, rotate=180] at (5,10) {pentominoL={-1}};
         \pic[scale=\scfa,rotate=90] at (8,1) {pentominoT};
         \pic[scale=\scfa] at (8,7) {pentominoF};
        \pic[scale=\scfa] at (7,0) {pentominoX};
        \pic[scale=\scfa,rotate=180] at (14,10) {pentominoV};
        \pic[scale=\scfa] at (16,6) {pentominoN={-1}};
        \pic[scale=\scfa,rotate=180] at (12,3) {pentominoP};
        \pic[scale=\scfa,rotate=270] at (12,3) {pentominoU};
        \pic[scale=\scfa] at (14,0) {pentominoW};
        \pic[scale=\scfa] at (16,2) {pentominoY={-1}};
        \pic[scale=\scfa] at (5,4) {pentominoZ};
     \end{tikzpicture}
    \caption[Third move: Camille moves the ``L'' tile. ]{Third move: Camille moves the ``\,\begin{tikzpicture}[baseline={(current bounding box.center)}]
\pic[scale=0.1,rotate=0] at (0,0.5) {pentominoL};
\end{tikzpicture}'' tile and flips it. The fence now encloses 47 tiles. 21 moves remain.}
    \end{subfigure}
    \caption{An example of a game between Andréa, Billie and Camille. The view is zoomed in on a $14\times 14$ segment of the $20\times 20$ board.}\label{fig:game_example}
\end{figure}
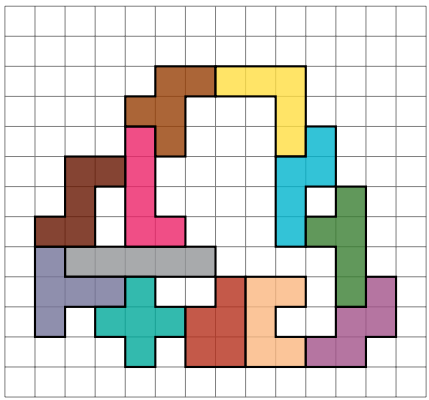
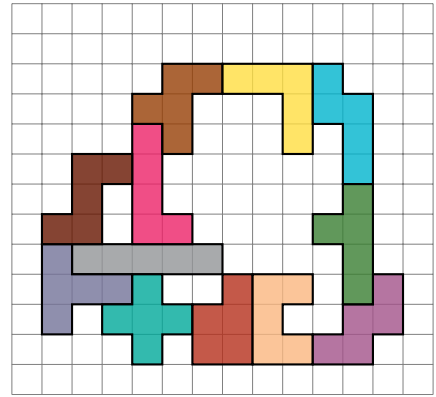
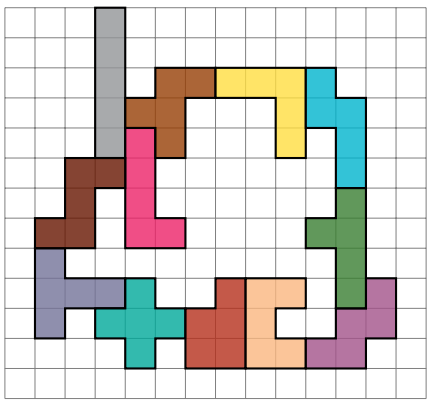
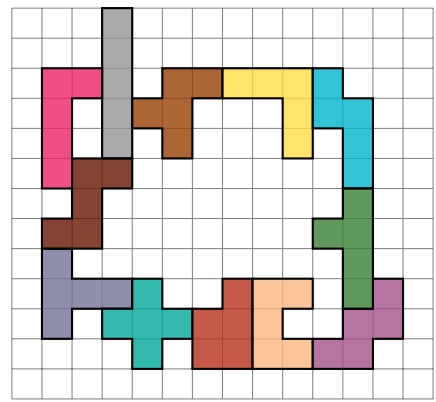

\pagebreak 

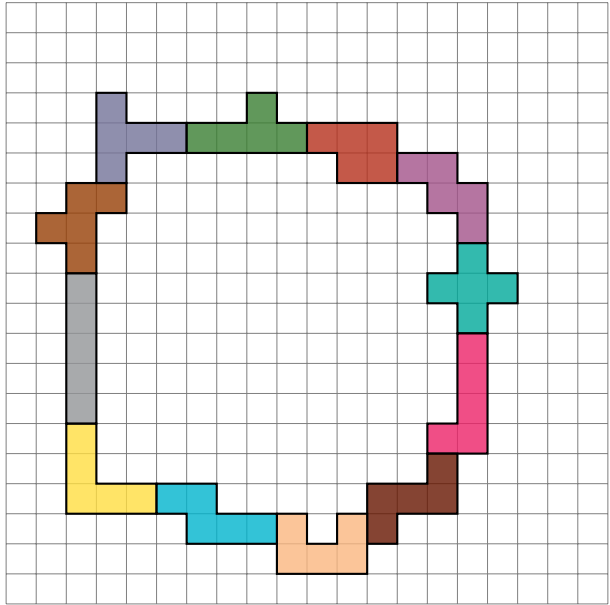
\begin{figure}[hp]
    \centering
    \begin{tikzpicture}[scale=\scfa]
    \draw[help lines] (0,0) grid (20,20);
        \pic[scale = \scfa,rotate=0] at (2,3){pentominoV};
        \pic[scale=\scfa,rotate=0] at (2,6) {pentominoI};
        \pic[scale = \scfa,rotate =0] at (1,11) {pentominoF};
        \pic[scale = \scfa,rotate=90] at (6,14){pentominoT};   \pic[scale = \scfa,rotate=270] at (6,15){pentominoY={-1}};
        \pic[scale = \scfa,rotate=270] at (10,16) {pentominoP};
        \pic[scale = \scfa,rotate=90] at (16,12) {pentominoW};
        \pic[scale = \scfa,rotate=0] at (14,9) {pentominoX};
        \pic[scale = \scfa,rotate=0] at (16,5) {pentominoL={-1}};
        \pic[scale = \scfa,rotate=270] at (12,2) {pentominoZ={-1}};
        \pic[scale = \scfa,rotate=0] at (9,1) {pentominoU};
        \pic[scale = \scfa,rotate=90] at (9,2) {pentominoN};
    
    \end{tikzpicture}
    \caption{Example of one of the optimal solutions of 128 tiles.}
    \label{fig:solmax}
\end{figure}
\end{document}